\documentclass{amsart}
\usepackage[T1]{fontenc}
\usepackage[latin1]{inputenc}
\usepackage{hyperref}
\makeatletter

\theoremstyle{plain}    
\newtheorem{thm}{Theorem}[section]
\numberwithin{equation}{section} 
\numberwithin{figure}{section} 
\theoremstyle{plain}    
\newtheorem{cor}[thm]{Corollary} 
\newtheorem{lem}[thm]{Lemma} 
\theoremstyle{plain}    
\newtheorem{prop}[thm]{Proposition} 
\theoremstyle{remark}
\newtheorem{rem}[thm]{Remark}
\theoremstyle{remark}    
\newtheorem{notation}[thm]{Notation} 

\usepackage[arrow,matrix,curve,ps]{xy}
\xyoption{dvips}
\CompileMatrices

\renewcommand{\H}{\mathcal H}

\renewcommand{\O}{\mbox{$\mathcal{O}$}}
\renewcommand{\P}{\mbox{$\mathbb{P}$}}

\DeclareMathOperator{\Aut}{Aut}
\DeclareMathOperator{\birat}{birat}
\DeclareMathOperator{\Chow}{Chow}

\DeclareMathOperator{\Hom}{Hom}
\DeclareMathOperator{\Image}{Image}
\DeclareMathOperator{\locus}{locus}

\DeclareMathOperator{\rank}{rank}
\DeclareMathOperator{\RatCurves}{RatCurves}

\DeclareMathOperator{\Span}{Span}

\DeclareMathOperator{\Univ}{Univ}

\makeatother

\begin{document}

\title{Lines on contact manifolds}

\date{\today }

\author{Stefan Kebekus}



\thanks{The author gratefully acknowledges support by a
  Forschungsstipendium of the Deutsche Forschungsgemeinschaft}

\address{Stefan Kebekus, Institut für Mathematik, Universität
  Bayreuth, 95440 Bayreuth, Germany}

\email{stefan.kebekus@uni-bayreuth.de}

\urladdr{http://btm8x5.mat.uni-bayreuth.de/\~{}kebekus}

\maketitle
\tableofcontents

\section{Introduction}

If $X$ is a complex projective manifold which carries a contact
structure, then the results of \cite{Dem00} and \cite{KPSW00} show
that $X$ is either isomorphic to a projectivized tangent bundle of a
complex manifold, or that $X$ is Fano and $b_2(X)=1$. In this paper we
study the latter case where $X$ is Fano. It is generally believed that
these assumptions imply that $X$ is homogeneous ---for an
introduction, see the excellent survey in \cite{Bea99}.

It follows from our previous work \cite{Keb00a} that $X$ can always be
covered by lines. Thus, it seems natural to consider the geometry of
lines in greater detail. We will show that if $x\in X$ is a general
point, then all lines through $x$ are smooth. If $X\not \cong
\P_{2n+1}$, then the tangent spaces to these lines generate the
contact distribution at $x$. It follows that the contact structure on
$X$ is unique, thus answering a question of C.~LeBrun
\cite[question~11.3]{Leb97}. The result was previously obtained by
C.~LeBrun \cite{Leb95} if $X$ is a twistor space.

\subsubsection*{Acknowledgement }

The paper was written while the author enjoyed the hospitality of RIMS
in Kyoto.  The author is grateful to Y.~Miyaoka for the invitation and
to the members of that institute for creating a stimulating
atmosphere. He would like to thank S.~Helmke, J.-M.~Hwang, S.~Kovács
and Y.~Miyaoka for a number of enlightening discussions on the
subject. Finally, the author thanks the referee for criticism
concerning the exposition of earlier versions of the paper.

\section{Setup}
\label{sec:setup}

\subsection{Contact Manifolds}

Throughout the present work, let $X$ be a complex projective contact
manifold of dimension $2n+1$. By definition, the contact structure is
given by a sequence of vector bundles 
$$ 
\xymatrix { 0
\ar[r] & F \ar[r] & {T_X} \ar[r]^{\theta} & L \ar[r] & 0} 
$$
where $L$ is of rank 1 and the contact form $\theta \in
H^0(\Omega_X^1\otimes L)$ yields a nowhere vanishing section 
$$
\theta \wedge (d \theta)^{\wedge n}\in H^0(K_X\otimes L^{n+1})
$$
---it is an elementary calculation to see that $\theta \wedge (d
\theta)^{\wedge k}$ is a well-defined section of
$\Omega_X^{2k+1}\otimes L^{k+1}$ for all numbers $n\geq k\geq 1$. In
particular, we assume that $-K_X=(n+1)L$.

The assumptions imply that the natural map $[,]:F\otimes F\to L$
derived from the Lie-bracket is non-degenerate. In view of the
Frobenius theorem this means that if $Y\subset X$ is an $F$-integral
submanifold, i.e.~one where $T_Y \subseteq F|_Y$, then $\dim Y\leq n$.
If $Y$ is of maximal dimension $\dim Y=n$, then $Y$ is called
``Legendrian''. Note that some authors (e.g.~\cite{Hwa97},
\cite{KPSW00}) prefer to use ``Lagrangian'' instead of
``Legendrian''.

The usual Darboux theorem of real contact and symplectic geometry
applies equally well in the complex case. Thus, for any point $x\in X$
we can find coordinates $(z_i)_{i=1\ldots 2n+1}$ centered about $x$
and a bundle coordinate for $L$ so that we can write
$$
\theta = dz_{2n+1}+\sum_{i=1\ldots n} z_idz_{n+i}
$$
In particular, we remark the following:
\begin{rem}
  \label{rem:transversality}
  If $x \in X$ is any point and $\vec v \in T_{X,x}$ any tangent
  vector, then there exists a Legendrian submanifold $U\subset X$
  which contains $x$ and is transversal to $\vec v$, i.e.~$\vec v \not
  \in T_U|x$.
\end{rem}

\subsection{Parameter spaces}

For the benefit of readers coming from differential geometry we will
briefly recall some facts about the parameter spaces which we will use
in the sequel. Our chief reference will be \cite[chap.~II]{K96}, and
our notation will be compatible with this book. Mori's paper
\cite{Mor79} on the Hartshorne conjecture is also recommended for
these matters.

If $V$ is any projective manifold and $C$ a projective variety, then
we will often parameterize those morphisms from $C$ to $V$ which are
birational onto their images. In fact, there exists a scheme
$\Hom_{bir} (C,V)$ whose geometric points correspond to these
morphisms. Furthermore, there exists a ``universal morphism'': $\mu:
\Hom_{bir} (C,V)\times C \to V$. We refer to \cite[chapt.~II.1]{K96}
for an authoritative reference on this.

The tangent space to $\Hom_{bir} (C,V)$ is described as follows: If a
birational morphism $f:C\to V$ is given, then the (Zariski-)tangent
space of $\Hom_{bir}(C,V)$ at $f$ corresponds naturally to sections
in $H^0(C, f^*(T_V))$.

If $c\in C$ and $v \in V$ are (geometric) points, the subfamily of
morphisms mapping $c$ to $v$ is usually denoted as $\Hom_{bir}(C,
V,c\mapsto v)$. The tangent space to $\Hom_{bir}(C, V, c\mapsto v)$
at a point $f$ corresponds to $H^0(C, f^*(T_V)\otimes \mathcal J_c)$,
where $\mathcal J_c$ is the ideal sheaf of $c$.

In the special case that $C\cong \P_1$, and $f\in \Hom_{bir} (\P_1,
V)$ a morphism whose image contains a point $v \in X$, the Riemann-Roch
theorem yields an estimate for the dimensions of the deformation
spaces
\begin{equation}
  \label{eq:dim_estimate}
  \begin{array}{rcl}
    \dim_{[f]} \Hom_{bir} (\P_1,V) & \geq & -K_V.\ell + \dim V \\
    \dim_{[f]} \Hom_{bir} (\P_1,V,[0:1]\mapsto v) & \geq & -K_X.\ell
  \end{array}
\end{equation}
where $\ell := \text{Image}(f)$. See \cite[prop.~II.1.13 and
thm.~II.1.7]{K96} for an explanation and a proof.

The group $\P SL_2$ acts on the normalization $\Hom^n_{bir}(\P_1, V)$
of $\Hom_{bir}(\P_1, V)$, and the geometric quotient in the sense of
Mumford \cite{FM82} exists, see \cite[lem.~9]{Mor79}. More precisely,
by \cite[thm.~II.2.15]{K96} there exists a commutative diagram
\begin{equation}
  \label{diag:rat_curves}
  \xymatrix{ 
    \Hom^n_{bir} (\P_1, V)\times \P_1 \ar[d] \ar[r]^(.6){U} 
    \ar@/^.6cm/[rr]^{\mu}
    & {\Univ^{rc}(V)} \ar[r]^(.6){\iota} \ar[d]_{\pi} & {X} \\ 
    \Hom^n_{bir} (\P_1, V) \ar[r]^{u} & {\RatCurves^n(V)} }  
\end{equation}
where $u$ and $U$ are principal $\P SL_2$ bundles, $\pi$ is a
$\P_1$-bundle and the restriction of $\iota$ to any fiber of $\pi$ is
a morphism which is birational onto its image. We call the quotient
space $\RatCurves^n(V)$ the ``parameter space of rational curves on
$V$''.  The letter ``$n$'' in $\RatCurves^n$ may be a little
confusing. It has nothing to do with the dimension of $V$, but serves
as a reminder that the parameter space is isomorphic the {\em
  n}ormalization of a suitable quasiprojective subset of the
Chow-variety.

\subsection{Lines}
\label{sect:lines}

Unless otherwise mentioned, throughout this work we will assume that
$X$ is Fano and that $b_2(X)=1$. In this setup it follows from the
classic work of Mori (\cite{Mor79}, but see
also~\cite[lect.~1]{CKM88}) that we can find an irreducible component
$H \subset \RatCurves^n(X)$ with the following properties:
\begin{enumerate}
\item the evaluation morphism $\iota$ is dominant
  
\item if $x \in X$ is a general point, then the subfamily $H_x :=
  \pi(\iota^{-1}(x)) \subset H$ (i.e.~the subfamily which
  parameterizes curves containing $x$) is compact.
  
\item if $\ell \subset X$ is a curve which is associated with a point
  of $H$, then $1 \leq -K_X.\ell \leq \dim X+1$.
\end{enumerate}
Since $X$ is a contact manifold, we have that $-K_X=(n+1)L$, and it
follows from point (3) that either $L.\ell =2$ or $L.\ell =1$.  

If $L.\ell =2$, then the estimate~(\ref{eq:dim_estimate}) shows that
$$
\dim \Hom_{bir}(\P_1,X,[0:1]\mapsto x) \geq -K_X.\ell = \dim X+1.
$$
Because there is a 2-dimensional group of automorphisms of $\P_1$
which fix $[0:1]$, the assumption that $L.\ell = 2$ implies
$$
\dim H_x \geq \dim X-1,
$$
and it follows from the generalized Kobayashi-Ochiai theorem of
\cite[thm.~3.6]{Keb00a} or \cite[thm.~0.2]{KS99} that $X\cong
\P_{2n+1}$. For that reason we assume in the sequel that $L.\ell =1$.
We call $\ell$ a ``contact line''.

\begin{rem}
  \label{rem:compactness}
  The assumption that $L.\ell = 1$ implies that the irreducible
  component $H \subset \RatCurves^n(X)$ is compact. See
  \cite[prop.~II.2.14]{K96}.
\end{rem}

\begin{rem}
  \label{rem:lines_are_intl}
  Let $f:\tilde{\ell }\to \ell $ be the normalization. Since
  $\tilde{\ell }\cong \P_1$, and $T_{\tilde{\ell }}\cong
  \O_{\P_1}(2)$, it is clear that the map $T_{\tilde{\ell}} \to
  f^*(L)$ must be trivial. It follows that $T_{\ell}|_x\subset F|_x$
  for all smooth points $x\in \ell $. We say that contact lines are
  $F$-integral where they are smooth.
\end{rem}

\section{Deformations of lines}

We will show that lines passing through a general point are smooth.
For this, we employ deformations of lines in order to obtain sections
of $L$. The following local proposition shows that there are severe
restrictions for such sections to exist.

\begin{prop}
  \label{Prop:Defo_of_F_intl_curves}
  Consider a family of morphisms $\Phi_t:\Delta_C\to X$ written as
  $$
  \begin{matrix}
    \Phi : &\Delta _{\H } \times \Delta_C & \to      & X \\
           &  (t,z)                       & \mapsto  & \Phi_t(z)
  \end{matrix}
  $$
  where $\Delta_{\H}$ and $\Delta_C$ are unit disks. Assume that
  for all $t\in \Delta_{\H}$ the image $\Phi_t(\Delta_C)$ is
  $F$-integral, i.e.~assume that $\Phi^*(\theta )\left(\frac{\partial
      }{\partial z}\right)\equiv 0$. If
  $$
  \sigma =\Phi^*_0(\theta )\left( \frac{\partial }{\partial t}\right)
  \in H^0(\Delta_C,\Phi_0^*(L))
  $$
  vanishes at $0\in \Delta_C$ but does not vanish identically, then
  $\sigma $ vanishes with multiplicity at least two if and only if
  \begin{equation}
    \label{orthogonal}
    \Phi ^{*}(d\theta )\left( \left.\frac{\partial }{\partial
        t}\right|_{(0,0)},\left. \frac{\partial }{\partial z} 
    \right|_{(0,0)}\right) =0.
  \end{equation}
\end{prop}

\begin{rem}
  The exterior derivative $d\theta$ which appears in
  equation~(\ref{orthogonal}) depends on a choice of bundle
  coordinates for $L$ and is therefore not well-defined. Note,
  however, that the requirement $\sigma (0)=0$ implies that $T\Phi
  \left( \frac{\partial }{\partial t}|_{(0,0)}\right) \in F|_{\Phi
    (0,0)}$, where $T\Phi $ is the tangent map associated with $\Phi$.
  In this setting an elementary calculation shows that the validity of
  equation~(\ref{orthogonal}) does in fact not depend on the choice of
  bundle coordinates.
\end{rem}

\begin{proof}
  Recall the following formula: if $\omega $ is any 1-form on a
  manifold, and $\vec{X}_{0}$ and $\vec{X}_{1}$ are vector fields,
  then
  \begin{equation}
    \label{D_theta}
    d\omega (\vec{X}_{0},\vec{X}_1) =
    \vec{X}_0(\omega (\vec{X}_1)) - 
    \vec{X}_1(\omega (\vec{X}_2))- \omega ([\vec{X}_0,\vec{X}_1]).
  \end{equation}
  See e.g.~\cite[prop.~2.25(e) on p.~70]{War71} for an explanation. We
  choose local bundle coordinates on $\Phi ^{*}(L)$ and apply equation
  (\ref{D_theta}) with $\omega =\Phi ^{*}(\theta )$, $\vec{X}_0
  =\frac{\partial }{\partial z}$ and $\vec{X}_1 =\frac{\partial
    }{\partial t}$. Since $\frac{\partial }{\partial t}$ and
  $\frac{\partial }{\partial z}$ commute and $\Phi ^{*}(\theta )\left(
    \frac{\partial }{\partial z}\right) \equiv 0$, we obtain
  $$
  \frac{\partial }{\partial z}\Phi ^{*}(\theta )\left( \frac{\partial
      }{\partial t}\right) =
  d\Phi ^{*}(\theta )\left( \frac{\partial }{\partial
      t},\frac{\partial }{\partial z}\right) .
  $$
  Note that $d\Phi^*(\theta)=\Phi^*(d\theta)$ and evaluate at
  $t=0$, $z=0$.
\end{proof}

Another argument which uses the deformation of lines shows that most
lines are smooth.

\begin{prop}
  \label{Prop:Lines_through_genl_point_are_smooth}
  If $x\in X$ is a general point and $\ell$ a contact line passing
  through $x$, then $\ell$ is smooth.
\end{prop}

\begin{proof}
  Assume to the contrary, i.e.~assume that for a general point $x\in
  X$ there exists a singular line $\ell$ passing through $x$. Recall
  that the rational curve $\ell$ can always be dominated by an
  integral singular plane cubic, i.e.~by a rational curve with a
  single node or cusp. We will reach a contradiction by constructing a
  section of the pull-back of $L$ to the plane cubic which vanishes at
  a prescribed generically chosen point. This section will be
  constructed by a deformation of the singular curves.
  
  Because $x$ was chosen to be general, there exists a singular
  (i.e.~cuspidal or nodal) plane cubic $C\subset \P_2$ and an
  irreducible component $\H \subset \Hom^{\birat} (C,X)$ such that the
  universal morphism $\mu :\H \times C\to X$ is dominant and such that
  for all $f\in \H $ we have $\deg f^*(L)=1$.
  
  Fix a general morphism $f\in \H $ and note that there is an open set
  $U\subset C$ such that for all $c\in U$, the tangent map of the
  restricted morphism $\mu_c:=\mu |_{\H \times \{c\}}$ has maximal
  rank at $f$:
  $$
  \rank_{[f]} T\mu_c= \dim X = 2n+1.
  $$
  Recall from \cite[II.6.10.2, II.6.11.4 and Ex.~II.6.7]{Ha77} that
  the smooth points of $C$ are in 1:1-correspondence with line bundles
  of degree one, and fix a point $c\in U$ such that $\O_C(c)\not \cong
  f^*(L)$.
  
  Next, let $U_X \subset X$ be a Legendrian submanifold of $X$ which
  contains $x$ and is transversal to $f(C)$ at $x$. By
  remark~\ref{rem:transversality}, these exist in abundance. Since
  $\mu_c$ has maximal rank, we can find a section $U_{\H}\subset \H $
  over $U_X$, i.e.~a submanifold $U_{\H}$ such that
  $\mu_c|_{U_{\H}}:U_{\H }\to U_X$ is an isomorphism. By construction,
  $\mu (U_{\H}\times C)$ has dimension $n+1$ and cannot be Legendrian.
  It follows that there exists a unit disc $\Delta_{\H}\subset U_{\H}$
  with coordinate $t$ centered about $f=\{t=0\}$ such that
  $$
  \sigma :=f^{*}(\theta )\left( \left. \frac{\partial }{\partial t}
    \right|_{t=0}\right) 
  \in H^{0}(C,f^{*}(L))\setminus \{0\}.
  $$
  For this, recall that the tangent vector $\frac{\partial
    }{\partial t}|_{t=0}\in T_{\H}|_f$ is canonically
  identified with an element in $H^0(C,f^*(T_X))$.  By choice of
  $U_{\H }$, we have $\sigma (c)=0$. By choice of $c$, this is
  impossible, a contradiction.
\end{proof}

It is an immediate corollary that a tangent morphism exists.

\begin{cor}
  If $x\in X$ is a general point, then there exists a tangent morphism
  $\tau_x:H_x \to \P(F^*|_x)$ which sends a line $\ell$ to its tangent
  space $T_\ell|_x \subset F|_x$.
\end{cor}

The morphism $\tau_x$ was already studied in \cite{Hwa97}. It was,
however, not all clear at that time that $\tau_x$ really is a morphism
and not just a rational map. See \cite{Keb00a} for a weaker but more
general result.

Finally, we remark that deformations of lines through a general point
are unobstructed in a strong sense.

\begin{lem}\label{lem:splitting_type}
  If $x\in X$ is a general point and $\ell$ any line through $x$, then
  $$
  T_X|_\ell \cong \O_\ell(2)\oplus \O_\ell(1)^{\oplus n-1} \oplus
  \O_\ell^{\oplus n+1}.
  $$
\end{lem}
\begin{proof}
  It follows from the definition of the contact structure that $F\cong
  F^*\otimes L$. Since $L|_\ell \cong \O_\ell(1)$, and since vector
  bundles on $\P_1$ always decompose into sums of line bundles we may
  therefore write
  $$
  F|_\ell \cong \bigoplus_{i=1}^n \left(
    \O_\ell(a_i)\oplus\O_\ell(1-a_i) \right)
  $$
  where $a_i > 0$. Thus, the splitting of $F|_\ell$ has exactly $n$
  positive entries. It follows that the splitting of $T_X|_\ell$ has
  at most $n$ positive entries. By \cite[prop.~1.1]{KMM92},
  $T_X|_\ell$ is nef, and since $c_1(T_X|_\ell)=n+1$, the claim
  follows.
\end{proof}

We will apply lemma~\ref{lem:splitting_type} to study the locus of
lines through a general point. For this, fix a general point $x\in H$,
define the subfamily $H_x\subset H$ as in section~\ref{sect:lines} and
consider the restricted diagram associated to
diagram~\ref{diag:rat_curves} on page~\pageref{diag:rat_curves}.
$$
\xymatrix{ {\tilde{U_x}} \ar[r]^(.3){\tilde \iota_x}
  \ar[d]_{\tilde \pi_x} & {\locus(H_x) \subset X} \\
  {\tilde {H_x}} \ar[r]^(.45){\tilde \tau_x} & {\P(F^*|_x)} }
$$
Here $\tilde{H_x}$ is the normalization of $H_x$, $\tilde{U_x}$ the
normalization of the pull-back $\Univ^{rc}(X) \times_{\RatCurves^n(X)}
{\tilde{H_x}}$ and $\locus (H_x) = \iota (\pi^{-1}(H_x))$. Recall from
remark~\ref{rem:compactness} that $H$ and therefore $H_x$ are compact.
In particular, $\locus(H_x)$ is a proper subvariety of $X$.

It follows from \cite[thms.~II.3.11.5 and II.2.8]{K96} that $\tilde
H_x$ is smooth and $\tilde \pi_x$ a $\P_1$-bundle. As an immediate
corollary to the preceding lemma, we obtain that both $\tilde \iota_x$
and $\tilde \tau_x$ are immersive.

\begin{cor}
  If $x\in X$ is a general point, then 
  \begin{enumerate}
  \item the universal morphism $\tilde\iota_x:\tilde U_x\to
    \locus(H_x)\subset X$ is a birational immersion away from a
    section $\sigma_\infty$ which is contracted to $x$.
  
  \item the tangent map $\tilde \tau_x:\tilde H_x \to \P(F^*|_x)$ is
    also an immersion
  \end{enumerate}
\end{cor}

\begin{proof}
  The fact that $\tilde \iota_x$ and $\tilde \tau_x$ are immersive
  follows from \cite[props.~II.3.4 and II.3.10]{K96} and
  lemma~\ref{lem:splitting_type}. It follows from an argument of
  Miyaoka that $\tilde \iota_x$ is birational because all lines
  through $x$ are smooth. For this, see~\cite[prop.~V.3.7.5]{K96}.
\end{proof}

\section{Lines through a fixed point}

For a better understanding of contact Fano manifolds, the locus of
lines through a given point is of greatest interest. The following
proposition gives a first description. This result is contained
implicitly in \cite[sect.~2]{KPSW00} and we could have used the
results of that paper here, but we prefer to give a short and
self-contained proof in our context.

\begin{prop}
  \label{Prop:equidim}
  If $x\in X$ is any point, then $\locus (H_x)$ has dimension $n$ and
  is $F$-integral where it is smooth.
\end{prop}

\begin{proof}
  Since $-K_X.\ell =n+1$, it follows from the
  estimate~(\ref{eq:dim_estimate}) for the dimension of the parameter
  space that
  $$
  \dim \Hom_{bir} (\P_1,X,[0:1]\mapsto x)\geq n+2.
  $$
  By Mori's bend-and-break argument \cite[thm.~II.5.4]{K96}, for a
  given point $y\in X\setminus \{x\}$, there are at most finitely many
  lines containing both $x$ and $y$. It follows that
  $$
  \dim \locus (H_x)=\dim \Hom (\P_1,X,[0:1] \mapsto x)-
  \dim \Aut (\P_1,[0:1])\geq n.
  $$
  
  \subsection*{Claim} The subvariety $\locus(H_x)$ is $F$-integral
  where it is smooth.

  \subsection*{Application of the claim} It follows
  immediately from Frobenius' theorem and from the non-degeneracy of
  the contact distribution that $\dim \locus(H_x) \leq n$, and we are
  done.
  
  \subsection*{Proof of the claim} Let $y\in \locus(H_x)$ be a general
  (smooth) point, $\ell\in H_x$ a curve which contains $x$ and $y$ and
  is smooth at $y$. By general choice of $y$, such a curve can always
  be found. Let $f:\P_1\to \ell$ be a birational morphism with
  $f([0:1])=x$ and $f([1:1])=y$. If
  $$
  \H\subset \Hom^{\birat}(\P_1,X,[0:1]\mapsto x)
  $$
  is an irreducible component of the reduced Hom-scheme which
  contains $f$, then we have that
  \begin{equation}
    \label{eq:tg_span}
    T_{\locus(H_x)}|_y \subset T_{\ell}|_y + 
    \Image(T_{\H}|_f \to T_X|_y)
  \end{equation}
  where $T_{\H}|_f$ is identified with $H^0(\P_1, f^*(T_X) \otimes
  \O_{\P_1}(-[0:1]))$ and the map $T_{\H}|_f \to T_X|_y$ is an
  application of the tangent map $Tf$ and evaluation at $y$. In other
  words, the tangent space $T_{\locus(H_x)}|_y$ at $y$ is spanned by
  the tangent space to the curve $\ell$ and by sections of $f^*(T_X)$
  which vanish at $[0:1]$. We refer to \cite[prop.~II.3.4]{K96} for a
  proof of~(\ref{eq:tg_span}).
  
  In remark~\ref{rem:lines_are_intl} we have already seen that $T_\ell
  \subset F|_\ell$ so that it suffices to prove that
  $$
  \text{Image}(T_{\H}|_f \to T_X|_y) \subset F|_y
  $$
  In other words, we have to show that if $\Delta_\H \subset \H$ is
  any unit disk centered about $f$ with coordinate $t$, then the
  section $\sigma \in H^0(\P_1,f^*(T_X))$ associated with
  $\frac{\partial }{\partial t}|_{t=0}$ is contained in
  $H^0(\P_1,f^*(F))$.  For this, note that
  proposition~\ref{Prop:Defo_of_F_intl_curves} asserts that the
  section $f^*(\theta )(\sigma )\in H^0(\P_1,f^*(L))$ has a zero at
  $0\in \P_1$ whose order is at least two. But since $\deg f^*(L)=1$,
  this implies that $f^*(\theta )(\sigma)$ vanishes identically. In
  particular, $\sigma \in H^0(\P_1, f^*(F))$.
  
  This proves that $T_{\locus (H_x)}|_y\subset F|_y$. Since $y$ was
  chosen generically, $\locus(H_x)$ is $F$-integral where it is
  smooth, and we are done.
\end{proof}

Under the assumptions spelled out in section~\ref{sec:setup}, if $x\in
X$ is a general point, then the contact distribution $F|_x$ is
generated by the tangent spaces to lines through $x$. Hence, it is
canonically given. Before starting the proof, however, it is
convenient to introduce the following notation first.

\begin{notation}
  Consider the incidence variety
  $$
  V := \{(x',x'')\in X\times X \ |\  x''\in \locus(x') \} \subset X\times X.
  $$
  An elementary calculation shows that $V$ is a closed subvariety
  of $X\times X$. We call $V$ the ``universal locus of lines through
  points''.
  
  Let $\pi_1,\pi_2: V \to X$ are the projection morphisms. Then for
  every $x\in X$ we have that (set-theoretically) $\pi_2
  (\pi_1^{-1}(x)) =\locus (H_x)$. It may well happen that
  $\pi_1^{-1}(x)$ is not reduced for special points $x\in X$.
  
  If $Y\subset X$ is a subset, we shall write $V|_Y$ for
  $\pi_1^{-1}(Y)$.
\end{notation}

\begin{lem}
  \label{lem:defo_of_incidence}
  Let $V\subset X\times X$ be the universal locus of lines through
  points which we defined above.  Let $\Delta$ be a unit disk with
  coordinate $t$ and $\gamma:\Delta\to X$ an embedding. Then there
  exists an open set $V^0\subset V|_{\gamma(\Delta)}$ such that
  $\pi_2(V^0)$ is a submanifold of dimension
  $$
  \dim \pi_2 (V^0) = n+1.
  $$
  In particular, by Frobenius' theorem, $\pi_2(V^0)$ is not
  $F$-integral. 
\end{lem}
\begin{proof}
  We have already seen in proposition~\ref{Prop:equidim} that
  $\pi_1|_V$ is equidimensional of relative dimension $n$.  Thus, $V$
  is a well-defined family of algebraic cycles over $X$ in the sense
  of \cite[I.3.10]{K96} and the universal property of the Chow-variety
  yields a map $\phi :X\to \Chow (X)$ such that $V$ is the pull-back
  of the universal family over $\Chow(X)$. Because $\dim \locus
  (H_x)=n<\dim X$, it is clear that the image of $\phi $ is not a
  point. Use the assumption that $b_2(X)=1$ to obtain that $\phi$ is
  actually a finite morphism.  Because two reduced algebraic cycles
  are equal if and only if their supports agree, it follows that for a
  given point $x_0\in X$, there are at most finitely many points
  $(x_i)_{i=1\ldots k}\subset X$ such that
  $$
  \locus(H_{x_0}) = \locus(H_{x_i}).
  $$
  In particular, if $V^0\subset V|_{\gamma(\Delta)}$ is an open set
  such that $\pi_2|_{V^0}$ is an embedding, then $\pi_2(V^0)$ has
  dimension $n+1$. Hence the claim.
\end{proof}

With these preparations we can now start the proof of the main theorem
of this work.

\begin{thm}
  \label{thm:Generation_of_F}
  If $x\in X$ is a general point, then $F|_x$ is spanned by the image
  of the tangent map $\tau_x$.
\end{thm}

\begin{proof}
  Our argument involves an analysis of the deformations of
  $\locus(H_y)$ which are obtained by varying the base point $y$. We
  shall argue by contradiction and assume that the assertion of the
  proposition is wrong. With this assumption we will construct a
  family of morphisms $\P_1 \to X$ which contradicts
  proposition~\ref{Prop:Defo_of_F_intl_curves}, and we are done.
  
  We will now produce a map $\gamma$ to which
  lemma~\ref{lem:defo_of_incidence} can be applied. Assuming that the
  statement of the proposition is wrong, we can find an analytic open
  neighborhood $U=U(x)\subset X$ and a subbundle $F'\subset F|_U$ such
  that 
  \begin{enumerate}
  \item For all $y\in U$, the vector spaces $F'|_y$ and $\Span(\Image
    \tau_y)\subset F|_y$ are perpendicular with respect to the
    non-degenerate form $F\otimes F\to L$ which comes with the contact
    structure.
  \item All lines through $U$ are smooth.
  \end{enumerate}
  After shrinking $U$, if necessary, let $\vec{v}\in H^{0}(U,F')$ be a
  nowhere-vanishing vector field.  Thus, if $y\in U$ is any point and
  $\ell \ni y$ is any line through $y$, then
  \begin{equation}
    \label{form:perp_assert}
    T_{\ell }|_y\subset \vec{v}(y)^{\perp },
  \end{equation}
  where ``$\perp $'' again means: perpendicular with respect to the
  non-degenerate form on $F$. Let $\Delta $ be a unit disc with
  coordinate $t$ and $\gamma :\Delta \to X$ be an integral curve of
  $\vec{v}$ with $\gamma (0)=x$.
  
  Now let $\H \subset \Hom_{bir}(\P_1,X)$ be the family of morphisms
  parameterizing the curves associated with $H$. Set
  $$ 
  \H_{\Delta} := \{f\in \H\ |\ f([0:1])\in \gamma(\Delta)\}.
  $$
  If $\mu_\Delta : \H_\Delta \times \P_1 \to X$ is the universal
  morphism, then it follows by construction that
  $$
  \mu_\Delta(\H_\Delta\times \P_1) = 
  \pi_2(V|_{\gamma(\Delta)})\supset\pi_2(V^0).
  $$
  In particular, since $\pi_2(V^0)$ is not $F$-integral, for a
  general point $(f,p)\in \H_\Delta\times \P_1$ there exists a tangent
  vector $\vec w\in T_{\H_\Delta\times \P_1}|_{(f,p)}$ such that the
  image of the tangent map is not in $F$:
  $$
  T{\mu_\Delta}(\vec w) \not \in  F
  $$
  Decompose $\vec w = \vec w'+\vec w''$, where $\vec w\in
  T_{\P_1}|_p$ and $\vec w''\in T_{\H_\Delta}|_f$. Then, since
  $f(\P_1)$ is $F$-integral, it follows that $T{\mu_\Delta}(\vec
  w')\in F$ and therefore 
  \begin{equation}
    \label{eq:wprimeprime}
    T{\mu_\Delta}(\vec w'')\not \in F.    
  \end{equation}
  As a next step, choose an immersion 
  $$
  \begin{array}{rrcl}
    \beta : & \Delta & \to     & \H_\Delta \\
            & t      & \mapsto & \beta_t
  \end{array}
  $$
  such that $\beta_0=f$ and such that
  $$
  T\beta\left(\left. \frac{\partial}{\partial t}\right|_{t=0}\right) 
  = \vec w''
  $$
  In particular, if $\sigma\in H^0(\P_1, f^*(T_X))$ is the section
  associated with $\vec w'' = T\beta(\frac{\partial}{\partial
    t}|_{t=0})$, and $\sigma' := f^*(\theta)(\sigma)\in
  H^0(\P_1,f^*(L))$, then the following holds:
  \begin{enumerate}
  \item it follows from (\ref{eq:wprimeprime}) and from
    \cite[prop.~II.3.4]{K96} that $\sigma'$ is not identically zero.
    
  \item at $[0:1]\in \P_1$, the section $\sigma$ satisfies
    $\sigma([0:1])\in f^*(T_{\gamma(\Delta)})\subset f^*(F)$. In
    particular, $\sigma'([0:1]) = 0$.
    
  \item If $z$ is a local coordinate on $\P_1$ about $[0:1]$, then it
    follows from (\ref{form:perp_assert}) that
    $\frac{\partial}{\partial z}|_{[0:1]}\in f^*(F)$ and
    $\sigma([0:1]) \in f^*(F')$ are perpendicular with respect the the
    non-degenerate form.
  \end{enumerate}
  
  Items (2) and (3) ensure that we can apply
  proposition~\ref{Prop:Defo_of_F_intl_curves} to the family $\beta_t$.
  Since the section $\sigma'$ does not vanish completely, the
  proposition states that $\sigma'$ has a zero of order at least two
  at $[0:1]$. But $\sigma'$ is an element of $H^0(\P_1,f^*(L))$, and
  $f^*(L)$ is a line bundle of degree one. We have thus reached a
  contradiction, and the proof of theorem~\ref{thm:Generation_of_F} is
  therefore finished.
\end{proof}

It follows that there are only two types of contact manifolds whose
structure is not unique.

\begin{cor}
  If $X$ is any complex projective contact manifold with more than one
  contact structure, then either $X\cong \P_{2n+1}$ or $X\cong \P
  (T_Y)$ for a manifold $Y$ whose tangent bundle $T_Y$ has bundle
  automorphisms.
\end{cor}

We refer to \cite[sect.~2.6]{KPSW00} for a study of the different
contact structures on $\P(T_Y)$.

\begin{proof}
  By \cite[prop.~3.1]{KPSW00}, the canonical bundle $\omega _{X}$ is
  not nef. But then it has already been shown in
  \cite[thm.~1.1]{KPSW00} that $X$ is automatically of type $\P (T_Y)$
  if $b_2(X)>1$. We may therefore assume without loss of generality
  that $X$ is Fano and that $b_2(X)=1$.
  
  Let $x\in X$ be a general point, and $\ell \ni x$ a minimal rational
  curve through $x$. It follows from the classical argument of Mori
  that $-K_{X}.\ell \leq \dim X+1$. Since $X$ is a contact manifold,
  $-K_{X}=(n+1)L$ so that $L.\ell \in \{1,2\}$. If $L.\ell =2$, then
  $X\cong \P_{2n+1}$. If $L.\ell =1$, then
  theorem~\ref{thm:Generation_of_F} shows that the contact
  distribution $F$ is canonically defined, hence unique.
\end{proof}

\bibliographystyle{alpha}

\end{document}